\documentclass[preprint,12pt]{article}

 \usepackage{graphics}
 \usepackage{graphicx}
 \usepackage{epsfig}
 \usepackage[all]{xy}
 \usepackage{cite}

\usepackage{amssymb}
\usepackage{amsthm,amsmath}

\usepackage{latexsym, epsfig, psfrag,eepic,colordvi,bm}
\usepackage{graphics,color}
\usepackage{amsmath,amsfonts,amssymb,amscd}
\usepackage[all]{xy}
\usepackage{epstopdf}

\usepackage{amsthm,amsmath,amssymb}

\usepackage{graphicx}

\theoremstyle{plain}
\newtheorem{theorem}{Theorem}[section]
\newtheorem{lemma}[theorem]{Lemma}
\newtheorem{corollary}[theorem]{Corollary}
\newtheorem{proposition}[theorem]{Proposition}

\theoremstyle{definition}
\newtheorem{definition}[theorem]{Definition}
\newtheorem{example}[theorem]{Example}

\theoremstyle{remark}

\newcommand{\bai}{\hspace{6pt}}

\begin{document}

		

\title{Separation probabilities and analogues of a Zagier-Stanley formula}

\author{Ricky X. F. Chen\\
	\small Biocomplexity Institute and Initiative, University of Virginia \\[-0.8ex]
	\small 995 Research Park Blvd, Charlottesville, VA 22911, USA\\
	\small\tt chen.ricky1982@gmail.com
}

\date{}
\maketitle


\begin{abstract}
In this paper, we first obtain some analogues of a formula of Zagier (1995) and Stanley (2011).
For instance, we prove that 
the number of pairs of $n$-cycles whose product has $k$ cycles and has $m$ given elements contained in distinct cycles (or separated) is given by
$$
\frac{2 (n-1)! C_m(n+1,k)}{(n+m)(n+1-m)}
$$ 
when $n-k$ is even, where $C_m(n,k)$ is the number of permutations of $n$ elements
having $k$ cycles and separating $m$ given elements.
As consequences, we obtain the formulas for certain separation probabilities due to Du and 
Stanley, answering a call of Stanley for simple combinatorial proofs.
Furthermore, we obtain the expectation and variance of the number of fixed points in the product of two random
$n$-cycles.

  

\end{abstract}


\section{Introduction}

Let $\mathfrak{S}_n$ denote the symmetric group on $[n]=\{1,2,\ldots, n\}$.
We shall use the following two representations of a permutation $\pi\in \mathfrak{S}_n$:\\
\emph{two-line form:} the top line lists all elements in $[n]$, following the natural order.
The bottom line lists the corresponding images of elements on the top line, i.e.,
\begin{eqnarray*}
	\pi=\left(\begin{array}{ccccccc}
		1&2& 3&\cdots &n-2&{n-1}&n\\
		\pi(1)&\pi(2)&\pi(3)&\cdots &\pi({n-2}) &\pi({n-1})&\pi(n)
	\end{array}\right).
\end{eqnarray*}
\emph{cycle form:} regarding $\langle \pi\rangle$ as a cyclic group, we represent $\pi$ by its
collection of orbits (cycles). The number of cycles of $\pi$ is denoted by $C(\pi)$.
The set consisting of the lengths of these disjoint cycles is called the cycle-type of $\pi$, and 
denoted by $ct(\pi)$.
We can encode this set as an integer partition of $n$.
An integer partition $\lambda$ of $n$, denoted by $\lambda \vdash n$,
can be represented by a non-increasing integer sequence $\lambda=\lambda_1
\lambda_2\cdots$, where $\sum_i \lambda_i=n$, or as $1^{a_1}2^{a_2}\cdots n^{a_n}$, where
we have $a_i$ of part $i$ and $\sum_i i a_i =n$. A cycle of length $k$ is called a $k$-cycle
and a permutation with only cycles of length two is called a fixed point free involution.

Separation probabilities for products of permutations were studied in Bernardi, Du, Morales, and Stanley~\cite{bdms}, where
a special case is concerned with the probability of having the elements in $[m]$
contained in distinct cycles of the product of a uniformly chosen $n$-cycle
and a permutation chosen uniformly randomly from the set of permutations of cycle-type $\lambda$.
Explicit formulas for computing the separation probabilities for the cases of $\lambda$ being $n^1$ and
$2^k$ were obtained. For instance, when $\lambda=n^1$, the separation probability is given by
\begin{align*}
\begin{cases}
\frac{1}{m!}, & \mbox{if $n-m$
is odd},\\
\frac{1}{m!}+\frac{2}{(m-2)! (n+1-m)(n+m)}, & \mbox{otherwise}.
\end{cases}
\end{align*}
This case was also previously obtained by Du and Stanley~\cite{stan2}.
For a general $\lambda$, the separation probability is encoded in a certain coefficient of a generating function involving
symmetric functions, which is usually hard to extract.
An earlier work of Stanley~\cite{stan3} studied the case $m=2$ and $\lambda=n^1$ addressing a conjecture by Bona~\cite{bona},
and Stanley asked for combinatorial proofs for these results~\cite{stan2, stan3}.
A combinatorial proof for the case $m=2$ was given in Cori, Marcus and Schaeffer~\cite{cms}. 

In another line of studies, Zagier~\cite{zag} and Stanley~\cite{stan3} have independently showed
that the number of $n$-cycles $s$ such that the product $(1\bai 2\bai \cdots \bai n)\, s$ has $k$ cycles is $\frac{2}{n(n+1)} C(n+1,k)$,
where $C(n,k)$ stands for the signless Stirling number of the first kind, i.e., the number of permutations on $[n]$ with $k$ cycles.
See also combinatorial proofs in~\cite{cms,fv,chr-1}.

In this paper, we first enumerate the pairs of $n$-cycles whose product has $k$ cycles and has the elements in $[m]$ separated (resp.~fixed) and obtain analogues of the above Zagier-Stanley formula.
Specifically, we show that these numbers are respectively given by 
$$
\frac{2 (n-1)! C_m(n+1,k)}{(n+m)(n+1-m)}, \quad \frac{2 (n-1)! \widehat{C}_m(n+1,k)}{(n-m)(n+1-m)},
$$ 
when $n-k$ is even, where $C_m(n,k)$ (resp.~$\widehat{C}_m(n,k)$) is the number of permutations on $[n]$
with $k$ cycles and the elements in $[m]$ being separated (resp.~fixed), i.e., analogues of $C(n,k)$.
When $m=0$, these results obviously reduce to the Zagier-Stanley formula.
Our approach is purely combinatorial and based on extending the plane permutation framework as discussed in Chen and Reidys~\cite{chr-1} in order 
to study hypermaps and genome rearrangement problems. 

As consequences, we are able to prove the formulas of the separation probabilities for $\lambda=n^1$ and a general $m$, which is probably the most 
simple combinatorial proof. We also prove that the isolation probability, i.e., the elements in $[m]$ are fixed points, is given by
$\frac{1}{m!} {n-1 \choose m}^{-1}$.
As consequences of the latter, we obtain that the expected number of fixed points
in the product of two uniformly random $n$-cycles has an elegant expression $\frac{n}{n-1}$
and we obtain the probability for the product to be fixed point free,
and combined with the symmetry property in~\cite{bdms} we obtain an elegant
formula for another kind of separation probabilities (which will be made precise later).

The outline of the paper is as follows. In Section~\ref{sec2}, we briefly review and extend the plane permutation framework. In Section~\ref{sec3}, we obtain the analogues of the Zagier-Stanley formula by relating our enumeration problem to 
a simple problem of counting exceedances of certain permutations. In Section~\ref{sec4}, we derive the separation and isolation probabilities of $m$ elements.
In Section~\ref{sec6}, we conclude the paper with some remarks.

\section{General formulas}\label{sec2}

The plane permutation framework has proven to be effective in studying hypermaps, graph embeddings and genome
rearrangement distances~\cite{chr-1, chr-2}.
In this section, we first review some notation and results about plane permutations from~\cite{chr-1}.
Then, we discuss our strategy for computing separation probabilities of $m$ elements by extending these results.

\begin{definition}
	A \emph{plane permutation} on $[n]$ is a pair $\mathfrak{p}=(s,\pi)$ where $s=(s_i)_{i=0}^{n-1}$
	is an $n$-cycle and $\pi$ is an arbitrary permutation on $[n]$.
	Given $s=(s_0~s_1~\cdots ~s_{n-1})$,
a plane permutation $\mathfrak{p}=(s,\pi)$ is represented by a two-row array:
\begin{equation}
\mathfrak{p}=\left(\begin{array}{ccccc}
s_0&s_1&\cdots &s_{n-2}&s_{n-1}\\
\pi(s_0)&\pi(s_1)&\cdots &\pi(s_{n-2}) &\pi(s_{n-1})
\end{array}\right).
\end{equation}
	The permutation $D_{\mathfrak{p}}$ induced by the diagonal-pairs (cyclically) in the array, i.e., for $0<i< n$,
	$D_{\mathfrak{p}}(\pi(s_{i-1}))=s_i$ and
$D_{\mathfrak{p}}(\pi(s_{n-1}))=s_0$, is called the \emph{diagonal} of $\mathfrak{p}$.
\end{definition}\label{2def1}

We sometime refer to $s,\, \pi, \, D_{\mathfrak{p}}$ respectively as the upper horizontal, the vertical and the diagonal.
Obviously, we have 
$D_{\mathfrak{p}}=
s \pi^{-1}$. 
It should be pointed out that, although as a cyclic permutation, there is no absolute left-right order for the elements in $s$,
in this paper, we generally assume there is a left-right order, with the leftmost element being $s_0$.

In a permutation $\pi$ on $[n]$, $i$ is called an \emph{exceedance} if $i<\pi(i)$ following the natural order and an \emph{anti-exceedance} otherwise. Exceedances and anti-exceedances are among
the most well-known permutation statistics.
Note that $s$ induces a linear order $<_s$,
where $a<_s b$ if $a$ appears before $b$ in $s$ from left to right (with the leftmost element $s_0$).
Without loss of generality, we always assume $s_0=1$ unless explicitly stated otherwise.
Then, these concepts can be generalized for plane permutations as follows:
\begin{definition}\label{2def2}
	For a plane permutation $\mathfrak{p}=(s,\pi)$, an element $s_i$ is called an
	\emph{exceedance} of $\mathfrak{p}$ if $s_i<_s \pi(s_i)$, and an \emph{anti-exceedance} if $s_i\ge_s \pi(s_i)$.
\end{definition}

In the following,
any comparison of elements in $s,~\pi$ and $D_{\mathfrak{p}}$ references the linear order $<_s$.
Obviously, each $\pi$-cycle contains at least one anti-exceedance as it contains
a minimum, $s_i$, for which $\pi^{-1}(s_i)$ is an anti-exceedance. We call these trivial anti-exceedances
and refer to a \emph{non-trivial anti-exceedance} as an NTAE. Furthermore, in any cycle of length
greater than one, its minimum is always an exceedance.

\begin{example}
	For the plane permutation 
	\begin{equation*}
	\mathfrak{p}=\left(\begin{array}{cccccc}
	1& 3& 6 &2 & 5 & 4\\
	5& 4 &1&3 & 6 & 2
	\end{array}\right),
	\end{equation*}
	$3$ is an exceedance, $1$ is a trivial anti-exceedance, and $5$ is an NTAE.
\end{example}

Let $\mathfrak{p}=(s,\pi)$ be a plane permutation. 
A \emph{diagonal block} of $\mathfrak{p}$ is a set of consecutive diagonal-pairs.
A \emph{transposition action} on the diagonal of $\mathfrak{p}$
transposes two adjacent diagonal blocks of $\mathfrak{p}$.
Specifically,
for a sequence $h=(i,j,k)$ such that $i\leq j<k$
and $\{i,j,k\}\subset [n-1]$, if we transpose the two diagonal-blocks determined by the continuous segments $[s_i,s_j]$ and $[s_{j+1},s_k]$, we obtain a new two-row array $\mathfrak{p}^h=(s^h,\pi^h)$:
\begin{eqnarray*}
\left(
\vcenter{\xymatrix@C=0pc@R=1pc{
\cdots & s_{i-1}  & s_{j+1}\ar@{--}[dl] &\cdots & s_{k-1}&\bai s_k\ar@{--}[dl] \bai\bai& s_i\ar@{-}[dl] &\cdots & s_{j-1} & s_{j}\ar@{-}[dl] & s_{k+1}  &\cdots\\
\cdots  & \pi(s_{j}) & \pi(s_{j+1}) & \cdots & \pi(s_{k-1})& \pi(s_{i-1}) & \pi(s_i) & \cdots\hspace{-0.5ex} & \pi(s_{j-1}) & \pi(s_{k})  & \pi(s_{k+1})&\cdots
}}
\right).
\end{eqnarray*}
Comparing $\mathfrak{p}$ and $\mathfrak{p}^h$, we have the following observations:
\begin{itemize}
	\item they have the same diagonal;
	\item the upper horizontals $s$ and $s^h$ differ by a transposition of the two continuous segments $[s_i,s_j]$ and $[s_{j+1},s_k]$;
	\item the maps $\pi$ and $\pi^h$ only differ at the images of the elements $s_{i-1}$, $s_j$, and $s_k$.
\end{itemize}
Thus, the transposition actions on the diagonal provide
a natural viewpoint on how different factorizations of the diagonal
into a long cycle (the upper horizontal) and another permutation (the vertical)
relate to each other. 
In particular,
the above last bullet implies that all components other than those containing the mentioned three elements of $\pi$ will be completely carried over to $\pi^h$ without any changes.
For those components containing the three elements, the three elements serve as certain breakpoints, where
the induced segments will be re-pasted in a certain way, depending on the distribution of the elements $s_{i-1}$, $s_j$, and $s_k$ in the components of $\pi$. 
Note that $\pi$ and $\pi^h$ must have the same parity.
Thus, the difference of the number of cycles in $\pi^h$ and $\pi$ is contained in $\{2, 0,-2\}$.
The NTAEs of $\mathfrak{p}$ can help us
to identify the transposition actions which change the number of cycles
by exactly two.

Let $D$ be a fixed permutation on $[n]$.
We consider the number of factorizations of $D$ into 
a long cycle $s$ on $[n]$ and a permutation $\pi$ on $[n]$ with
$k$ disjoint cycles in total, i.e.,~$D=s\pi$, which is related to enumerating one-face maps (when
$D$ is a fixed point free involution) and hypermaps.
Obviously, it is equivalent to considering the set of plane permutations
$\mathfrak{p}=(s, \pi^{-1})$ such that the diagonal is $D$
and the vertical has $k$ cycles.
Denote this set by $\tilde{U}_{k}^D$. 
The following result has been obtained.

 \begin{proposition}[Chen\&Reidys~\cite{chr-1}]\label{prop:general}
 Let $\tilde{Y}_1$ be the set of pairs $(\mathfrak{p},\epsilon)$ where $\mathfrak{p}\in \tilde{U}_{k}^D$ and $\epsilon$
 is an NTAE of $\mathfrak{p}$. Let $\tilde{Y}_2$ be the set of plane permutations $\mathfrak{p}\in \bigcup_{j\geq 1}  \tilde{U}_{k+2j}^D$ where there are $2j+1$ marked cycles in $\mathfrak{p}$ if $\mathfrak{p}\in  \tilde{U}_{k+2j}^D$. Then there is a bijection between $\tilde{Y}_1$ and $\tilde{Y}_2$.
 \end{proposition}
 
  The main idea behind the above bijection can be briefly summarized here.
 From a given pair $(\mathfrak{p},\epsilon)$ in $\tilde{Y}_1$, the NTAE $\epsilon$ determines a transposition on the diagonal 
 of $\mathfrak{p}$ such that the vertical of the resulting plane permutation after the transposition is
 obtained by splitting the cycle containing $\epsilon$ of $\mathfrak{p}$ into three cycles.
Obviously, $\epsilon$ is still contained in one of the three cycles.
Depending on whether $\epsilon$ is still an NTAE of the resulting plane permutations, additional
transpositions can be applied until $\epsilon$ is not an NTAE anymore.
Eventually, the original cycle containing $\epsilon$ will split into $2j+1$ cycles for some $j>0$ which will be marked.
Conversely, from a given element in $\tilde{Y}_2$, there is a unique way to merge the marked cycles
into one single cycle and create an NTAE.
We refer to~\cite{chr-1} for details.

 We remark that the bijection was motivated by the vertex slicing/gluing bijection on one-face maps in Chapuy~\cite{chapuy}.
 However, once we had the two-row array formulation of plane permutations, 
 it was in fact the natural transposition action on the diagonal, or an even broader perspective, rearrangement of the diagonal-pairs, that were first studied, due to their clear potential applications to the block-interchange and reversal distances of genome sequences.
 It turned out that the slicing/gluing operations are hiding there as two particular cases among others (see, e.g.,~\cite[Lemma~$7$]{chr-1}),
 somehow resolving the mystery of the slicing/gluing bijection~\cite{chapuy}.

 From the discussed bijection above, new bijections can be derived if there are some appropriate
 restrictions on the set of plane permutations. The main results of this paper 
 are based on such derived bijections. Let $D$ be a fixed permutation on $[n]$,
 and let ${U}_{m,k}^D$ be the set of plane permutations
$\mathfrak{p}=(s, \pi)$ with $D$ being the diagonal where there are $k$ cycles in 
the vertical and the elements in $[m]$ are contained in distinct cycles of the vertical. 
By similar arguments as Proposition~\ref{prop:general}, we obtain the following proposition.

 \begin{proposition}\label{prop:restricted}
 Let ${Y}_1$ be the set of pairs $(\mathfrak{p},\epsilon)$ where $\mathfrak{p}\in {U}_{m,k}^D$ and $\epsilon$
 is an NTAE of $\mathfrak{p}$. Let ${Y}_2$ be the set of plane permutations $\mathfrak{p}\in \bigcup_{j\geq 1}  {U}_{m,k+2j}^D$ where there are $2j+1$ marked cycles in $\mathfrak{p}$ if $\mathfrak{p}\in  {U}_{m,k+2j}^D$ and among the marked cycles at most 
 one of them contains an element in $[m]$. Then there is a bijection between ${Y}_1$ and ${Y}_2$.
 \end{proposition}
 \proof 
 We merely point out that if we start with a cycle containing an element from the set $[m]$, then it will generate
 $2j$ cycles not containing any elements from $[m]$ and one cycle containing exactly one element from $[m]$.
 Hence, for the merging operation in the converse procedure, among the chosen $2j+1$ cycles (to be marked) at most one of them can contain an
 element from $[m]$.  \qed

Let $p^D_{m,k}=|{U}_{m,k}^D|$, ${U}_{m,k}^{\lambda}=\bigcup_{ct(D)=\lambda} {U}_{m,k}^D$, and $p^{\lambda}_{m,k}=|{U}_{m,k}^{\lambda}|$. It is obvious that if $\mathfrak{p} \in {U}_{m,k}^{\lambda}$ has $a$ exceedances, then $Ne(\mathfrak{p})=n-a-k$ where $Ne(\mathfrak{p})$ denotes the number of NTAEs in $\mathfrak{p}$.

Then, as consequences of Proposition~\ref{prop:restricted}, we have
\begin{corollary} 
Suppose $D\in \mathfrak{S}_n$ is of cycle-type $\lambda\vdash n$, and let $p^{\lambda,n,a}_{m,k}$ denote the number of plane permutations
$\mathfrak{p}\in U^{\lambda}_{m,k}$ (on $[n]$) such that $\mathfrak{p}$
has $a$ exceedances.
Then, we have
	\begin{align}
	\sum_{\mathfrak{p}\in {U}_{m,k}^D} Ne(\mathfrak{p}) &= \sum_{j\geq 1} \left[m{k+2j-m\choose 2j}+{k+2j-m\choose 2j+1}\right] p^D_{m,k+2j} \;,\\
	\sum_{a\geq 0} (n-k-a) p^{\lambda,n,a}_{m,k} &= \sum_{j\geq 1} \left[m{k+2j-m\choose 2j}+{k+2j-m\choose 2j+1}\right] p^{\lambda}_{m,k+2j}\label{eq:U} \;.
	\end{align}
\end{corollary} 				

We further remark that these equations in the above corollary are inherently 
filtered out (or avoided) if following the map (and bicolored map~\cite{chapuy,walsh2}) perspective.
Because on the one hand, all plane permutations corresponding to maps (i.e.,~the diagonal being a fixed point free involution) have the same fixed number
(roughly speaking, $2g$ for $g$ being the genus)
of NTAEs, such equations never appear in the first place; On
the other hand, these equations do not really provide `valid'
recurrences from the enumeration perspective. (In order for obtaining valid recurrences, we have to apply a
sort of `reflection principle' to clear the parameter `$a$', as will be shown shortly.)
However, these `invalid' recurrences are actually 
the most valuable ingredients to make our approach work.

We denote $\mu\rhd_{k} \lambda$ if $\mu \vdash n$ can be obtained from $\lambda \vdash n$ by splitting one part into $k$ parts,
or equivalently, $\lambda$ from $\mu$ by merging $k$ parts into one part.
Let $\kappa_{\mu,\lambda}$ be the number of different ways of merging $k$ parts of $\mu$
in order to obtain $\lambda$ provided that $\mu \rhd_k \lambda$.
Note that we differentiate two parts of $\mu$ even if the two parts are of the same value.
For example, for $\mu=1^2 2^2$ and $\lambda=1^1 2^1 3^1$, we have
$\kappa_{\mu, \lambda}=4$.				
				
Let $\pi$ be a fixed permutation and $w^{\pi}_{\eta}$ be the number of 
distinct factorizations of $\pi$ into a long cycle and a permutation of cycle-type $\eta$. 
Let $V^{\pi}_{\eta}$ denote the corresponding set of plane permutations.

\begin{proposition}[Chen\&Reidys~\cite{chr-1}]\label{prop:V}
\begin{align}
	\sum_{\mathfrak{p}\in V^{\pi}_{\eta}} {Ne}(\mathfrak{p}) &=\sum_{j\geq 1} \sum_{\mu \rhd_{2j+1} \eta} \kappa_{\mu,\eta} w_{\mu}^{\pi} \; .
\end{align}
\end{proposition}

\begin{lemma}[Chen\&Reidys~\cite{chr-1}]\label{reflection}
Let $\mathfrak{p}=(s,\pi)$ be a plane permutation with diagonal $D_{\mathfrak{p}}$,
and let $\mathfrak{p'}=(s^{-1}, D_{\mathfrak{p}}^{-1})$. Then,
\begin{align}
Ne(\mathfrak{p})+Ne(\mathfrak{p'})=n+1-C(\pi)-C(D_{\mathfrak{p}}) \; .
\end{align}
\end{lemma}

For an integer partition $\lambda$, we denote the number of non-zero parts in $\lambda$ by $\ell(\lambda)$.
Now we can obtain a `valid' recurrence.
\begin{theorem}\label{thm:gen1}
For $\lambda \vdash n$ and $n+1-\ell(\lambda)-k>0$, we have
	\begin{align}
	 	 p^{\lambda}_{m,k}=\frac{\sum_{j\geq 1} \left[m{k+2j-m\choose 2j}+{k+2j-m\choose 2j+1}\right] p^{\lambda}_{m,k+2j}+\sum_{j\geq 1} \sum_{\mu \rhd_{2j+1} \lambda} \kappa_{\mu,\lambda} p^{\mu}_{m,k}}{n+1-\ell(\lambda)-k} \; .
	\end{align}
\end{theorem}
\proof Note that the set 
$$
V=\{(s^{-1}, D_{\mathfrak{p}}^{-1})\mid \mathfrak{p}=(s,\pi)\in U^{\lambda}_{m,k}\}=\bigcup_{\gamma}V^{\gamma}_{\lambda} \; ,
$$
where $\gamma$ is over all permutations with $k$ cycles such that the elements in $[m]$
are in distinct cycles.
According to Proposition~\ref{prop:V}, we have
\begin{align}\label{eq:V}
\sum_{\mathfrak{p}\in V} Ne(\mathfrak{p}) =\sum_{j\geq 1} \sum_{\mu \rhd_{2j+1} \lambda} \kappa_{\mu,\lambda} \sum_{\gamma}w_{\mu}^{\gamma}=\sum_{j\geq 1} \sum_{\mu \rhd_{2j+1} \lambda} \kappa_{\mu,\lambda} p^{\mu}_{m,k} \; .
\end{align}
Due to the one-to-one correspondence between $U^{\lambda}_{m,k}$ and $V$ as well as Lemma~\ref{reflection}, we have
$$
\sum_{\mathfrak{p}\in {U}_{m,k}^{\lambda}} Ne(\mathfrak{p}) + \sum_{\mathfrak{p}\in V} Ne(\mathfrak{p}) =(n+1-\ell(\lambda)-k) p^{\lambda}_{m,k} \; .
$$
Then, combining eq.~\eqref{eq:V} and eq.~\eqref{eq:U} completes the proof. \qed

It is known that in a plane permutation $\mathfrak{p}=(s,\pi)$, we have $C(\pi)+C(D_{\mathfrak{p}})\leq n+1$ where
the equality is attainabe.
See~\cite{chr-1, chr-2} for instance.
Thus, in order for obtaining an explicit formula for $p^{\lambda}_{m,k}$ for every $k$, it suffices 
to obtain $p^{\lambda}_{m, n+1-\ell(\lambda)}$ and the rest of the section is mainly devoted for this purpose.

It is observed in Walsh~\cite{walsh2} that there is some correspondence between maps and hypermaps.
In the following, we present a correspondence analogous to maps-hypermaps correspondence between plane permutations on $[n]$ and plane permutations on 
$[n]\bigcup [\bar{n}]$ where $[\bar{n}]=
\{\bar{1},\bar{2},\cdots, \bar{n}\}$. We remark that depending on the particular purposes, the construction between 
these two sets of plane permutations may be slightly different.

Let
\begin{equation*}
\mathfrak{p}=(s,\pi)=\left(\begin{array}{ccccc}
s_0& s_1 &\cdots &s_{n-1}& s_n\\
\pi(s_0)&\pi(s_1)&\cdots &\pi(s_{n-1}) &\pi(s_n)
\end{array}\right)
\end{equation*}			
be a plane permutation on $[n]$.
We associate a plane permutation $(\hat{s},\hat{\pi})$ on $[n]\bigcup [\bar{n}]$ where the diagonal is 
a fixed point free involution via the following procedure:
\begin{itemize}
	\item for any $i\in [n]$, put $\bar{i}$ right behind $i$ in the upper horizontal of $\mathfrak{p}$; and
	\item fill an appropriate element in $[\bar{n}]$ right below each $\bar{i}$ such that the diagonal of the resulting plane permutation is a fixed point free involution. 
\end{itemize}
It should not be hard to verify that there is a unique way to complete the second step above.
Thus, the associated plane permutation is unique.
The following is an example to illustrate this construction:
\begin{equation*}
\mathfrak{p}=\left(\begin{array}{cccccc}
1&2&3&4&5&6\\
4&5&6&1&2&3
\end{array}\right)\Longleftrightarrow
\left(\begin{array}{cccccccccccc}
1&\bar{1}&2&\bar{2}&3&\bar{3}&4&\bar{4}&5&\bar{5}&6&\bar{6}\\
4&\bar{5}&5&\bar{6}&6&\bar{1}&1&\bar{2}&2&\bar{3}&3&\bar{4}
\end{array}\right)=(\hat{s},\hat{\pi}).
\end{equation*}
It is obvious that the restriction of $\hat{\pi}$ to $[n]$
is the same as $\pi$. Regarding $\hat{\pi}$ on $[\bar{n}]$, we have
\begin{lemma}
	The restriction of $\hat{\pi}$ to $[\bar{n}]$,
	$\hat{\pi}|_{[\bar{n}]}$, has the same cycle-type as	$D_{\mathfrak{p}}$.	
\end{lemma}
\proof 
Note that identifying $\pi(i)$ and $\bar{i}$, $\hat{\pi}(\bar{i})$ and $s(i)$, i.e.,~the two elements on a same diagonal-pair, will preserve the cycle structure in the sense that
$$
\hat{\pi}:\bar{i}\rightarrow \hat{\pi}(\bar{i}), \quad D_{\mathfrak{p}}: \pi(i)\rightarrow s(i).
$$
Therefore, $\hat{\pi}|_{[\bar{n}]}$ has the same cycle-type as
$D_{\mathfrak{p}}$. 
\qed

For the example given above, we can check that $D_{\mathfrak{p}}=(153)(264)$ and
$\hat{\pi}|_{[\bar{5}]}=(\bar{1}\bar{5}\bar{3})(\bar{2}\bar{6}\bar{4})$ have the same cycle-type.
It is obvious that there are two types of cycles in $\hat{\pi}$, one has elements from $[n]$
and the other has elements from $[\bar{n}]$, and each diagonal-pair has exactly one element from $[n]$ and 
one from $[\bar{n}]$. 

Conversely, given a plane permutation $\hat{\mathfrak{p}}=(\hat{s},\hat{\pi})$ on $[n]\bigcup [\bar{n}]$
where $\hat{s}=( \cdots \bai i \bai \bar{i}\bai \cdots)$,
if the elements in any cycle of $\hat{\pi}$ are either from $[n]$ or from $[\bar{n}]$,
and the diagonal is a fixed point free involution,
then we can uniquely associate with it a plane permutation on $[n]$
where the diagonal has the same cycle type as $\hat{\pi}\mid_{[\bar{n}]}$
by just deleting the columns having elements from $[\bar{n}]$.

Based on the construction above, we can conclude

\begin{proposition}
Let $\hat{U}^{\lambda}_{m,\eta}$ be the set of plane permutations $(\hat{s},\hat{\pi})$ on $[n]\bigcup [\bar{n}]$
where the diagonal is a fixed point free involution with exactly one element from $[n]$
in each diagonal-pair, $\hat{s}$ is of the form $(\cdots \bai i \bai \bar{i} \bai \cdots)$,
$\hat{\pi}\mid_{[n]}$ is of cycle-type $\eta$ with the elements in $[m]$ in distinct cycles, and $\hat{\pi}\mid_{[\bar{n}]}$ is of cycle-type $\lambda$.
Let $\hat{U}^{\lambda}_{m,k}=\bigcup_{\eta,\, \ell(\eta)=k} \hat{U}^{\lambda}_{m,\eta}$.
Then, there is a one-to-one correspondence between $U^{\lambda}_{m,\eta}$ and $\hat{U}^{\lambda}_{m,\eta}$, as well as 
between $U^{\lambda}_{m,k}$ and $\hat{U}^{\lambda}_{m,k}$.
\end{proposition}

It is well known that when $\ell(\eta)=n+1-\ell(\lambda)$, the underlying structure of each element in $\hat{U}^{\lambda}_{m,\eta}$
corresponds to a tree. Thus, the set $\hat{U}^{\lambda}_{m,\eta}$ corresponds to a set of certain half-edge labelled plane trees, where by 
a half-edge we mean an end of an edge.

In a plane tree, we define the level of a vertex $v$ to be the length of the path from the root of the plane tree to the vertex $v$, and 
the root is on level $0$. The vertices on level $i+1$ adjacent to a vertex $v$ on level $i$ are called the 
children of $v$, and $v$ is called the parent of these said vertices. The degree of a vertex is the total number of edges incident to the vertex, while the outdegree of a vertex
is the number of children of the vertex, i.e., one less than the degree for a non-root vertex, and for the root, its degree is the same as its outdegree.
A non-root vertex of degree one is called a leaf. Non-leaf vertices are called internal vertices.
The following proposition describes the corresponding half-edge labelled plane trees.

\begin{proposition}
 For $\ell(\eta)=n+1-\ell(\lambda)$, the elements in $\hat{U}^{\lambda}_{m,\eta}$ are one-to-one corresponding to the set of plane trees $T$ of $n$ edges, where 
 the degree distribution of the vertices on the even level is $\eta$, the degree distribution of the vertices on the odd levels is $\lambda$, and the half-edges of $T$ 
 have labels constituting the set $[n]\bigcup [\bar{n}]$ in the following way: the half-edges incident to the even level vertices have labels in $[n]$ while the labels in $[m]$ are incident to distinct even level vertices, the half-edges incident to the odd level vertices have labels in $[\bar{n}]$, the rightmost half-edge incident to the root vertex has label $1$, and the half-edge paired with the counterclockwise neighbor of a half-edge with label $i$ is labelled by $\bar{i}$
 for any $i\in [n]$. 
 \end{proposition}
\proof 
For each plane permutation $(\hat{s},\hat{\pi}) \in \hat{U}^{\lambda}_{m,\eta}$,
we put the cycles of $\hat{\pi}$ as vertices and the elements in a cycle will 
be the half-edge labels (counterclockwisely) around the corresponding vertex, and we connect all half-edges according to the pairing 
relation specified 
by the diagonal $D$ of the given plane permutation. As a consequence, we obtain a graph $G$. Next, we 
can show that $G$ is a tree, i.e., a connect graph with $n+1$ vertices and $n$ edges.
Obviously, only connectedness needs to be verified.
Note that obtaining $D\hat{\pi}(x)$ is equivalent to, in $G$, starting with the half-edge $x$ (we
identify a half-edge and its label), and traveling to the counterclockwise 
neighbor of $x$, and going along the met edge to the other half-edge of the edge.
In this way, $\hat{s}=D\hat{\pi}$ being a long cycle on $[n]\bigcup [\bar{n}]$ implies that $G$
is connected whence being a tree. 
Next, if we view the vertex incident to half-edge $1$ as the root vertex of the tree such that half-edge $1$
is the rightmost incident half-edge, then the level of every vertex is uniquely determined and the left-to-right
relation among the vertices on the same level is uniquely determined as well. 
So the resulting structure is a half-edge labelled plane tree satisfying the condition 
specified in the proposition. It is clear how to reverse the above construction, whence the proposition.
 \qed
 
 It should be noted that by construction the labels of the half-edges incident to the odd level vertices are uniquely determined by those of the even level vertices. Thus, we could just ignore the labels of the half-edges incident to the odd level vertices in 
 the following.
 
 \begin{theorem}\label{thm:s-ini}
 Let
 $\ell(\mu)=d$, $\ell(\lambda)=t$ where $d+t=n+1$, and let $\ell_1^{\mu}=|\{\mu_i : \mu_i>1\}|$ and
 $\lambda=1^{a_1}2^{a_2} \cdots n^{a_n}$.
Then
  \begin{align}
{p}^{\lambda}_{m, \mu}=\frac{(t-1)! (d-1)! (n-m)!}{\prod_{i\geq 1} a_i !  (d-m)!}      \sum_{\substack{r>0, \; b\geq 0, \\ (r_1, \ldots, r_{b}), \\
(q_1,\ldots, q_{\ell_1^{\mu}-b+\bar{\delta}_{1r}})\\}}  {d-m \choose \ell_1^{\mu}-b-\bar{\delta}_{1r}}  {{m-1 \choose b} r \prod_{j\geq 1} (r_j+1) } \; ,
 \end{align}
  where $\bar{\delta}_{1r}=1$ if $r\neq 1$ and $0$ otherwise, and the tuples in the sum satisfy
 \begin{align*}
   \{r\} \bigcup \{r_1, \ldots, r_{b}\} \bigcup \{q_1,\ldots, q_{\ell_1^{\mu}-b-\bar{\delta}_{1r}}\} =\{r\} \bigcup \bigl(\{\mu_i-1 : \mu_i>1\} \setminus \{r-1\} \bigr),
 \end{align*}
 and $r=\mu_i$ for some $i$.
 \end{theorem}
\proof 
Let $\mathbb{T}$ be the set of vertex labelled plane trees,  where in a tree $T\in \mathbb{T}$, the labels for the even level vertices constitute $[d]$, the root always has label $d$, and the labels for the odd level 
vertices constitute $[\bar{t}]$, and the underlying structure of $T$ with the vertex labels ignored corresponds to an element in $\hat{U}^{\lambda}_{m,\mu}$.
Note that the labels for the half-edges are independent of the labels for the vertices.
Let $\hat{p}^{\lambda}_{m, \mu}=|\hat{U}^{\lambda}_{m,\mu}|$.
Then, $|\mathbb{T}|=(d-1)! t! \hat{p}^{\lambda}_{m,\mu}=(d-1)! t!  {p}^{\lambda}_{m, \mu}$.
In the following, we enumerate the number of trees in $\mathbb{T}$,
and we sometime refer to a vertex just by its label.
To facilitate the enumeration, we employ the following variation of Chen's algorithm~\cite{bill} on uniformly labelled plane trees.

For each $T\in \mathbb{T}$, we first decompose $T$ into a set of fibers (i.e.,~small labelled plane trees) according to the following procedure.
\begin{itemize}
\item[i.] Set $i_e=0,\, i_o=0$;
\item[ii.] We assume $\bar{j}<\overline{j+1}$ and any element $i\in [d]$ is smaller than $\bar{j}\in [\bar{t}]$.
In $T$, find the minmum internal vertex $v$ (in terms of its label) from $[d] \bigcup [\bar{t}]$ whose children are leaves.
Remove the fiber determined by $v$ (i.e.,~$v$ and its children) with all labels carried over, except for, in the case that $v$ is an even 
level vertex, remember the label of the half-edge incident to $v$ that is an end of the edge between $v$ and its parent;
\item[iii.] If $v$ has a label in $[d]$, place a vertex with label $(d+i_e+1)^*$ in the remaining tree at the original position of $v$ in $T$, and update $T$ as 
the resulting tree, and set $i_e=i_e+1$; If $v$ has a label in $[\bar{t}]$, place a vertex with label $(\overline{t+i_o+1})^*$ in the remaining tree at the original position of $v$ in $T$, and update $T$ as 
the resulting tree, and set $i_o=i_o+1$;
\item[iv.] If $T$ is not a fiber rooted on $d$, go to ii and continue, the procedure terminates otherwise.
\end{itemize}
In the end, we obtain a set of fibers, and we observe:
\begin{itemize}
\item In a fiber rooted on a vertex $v$ with a label from $[\bar{t}]$, the labels for the children of $v$ are from the set $[d-1]\bigcup 
\{ (d+1)^*, (d+2)^*,\ldots, (d+\ell_1^{\mu})^*\}$ if the root of $T$ has degree one or from the set $[d-1]\bigcup 
\{ (d+1)^*, (d+2)^*,\ldots, (d+\ell_1^{\mu}-1)^*\}$ if the root of $T$ has degree greater than one, and
the incident half-edges of the children having labels with the symbol $^*$ (starred labels) have labels the same as those `remembered'
ones.
The sizes of these fibers (i.e.,~the number of edges) rooted on the vertices with labels from $[\bar{t}]$
constitute the multiset $\{\lambda_i-1 : \lambda_i>1\}$. We shall refer to this subset of fibers as type-O fibers.
\item In a fiber rooted on a vertex $v$ with a label from $[d]$, the labels for the children of $v$ are from $[\bar{t}]\bigcup 
\{ (\overline{t+1})^*, (\overline{t+2})^*,\ldots, (\overline{t+\ell_1^{\lambda}})^*\}$, and besides a label from $[n]$ for each incident half-edge of $v$, 
there is one additional label from $n$ (i.e., the `remembered' one in step ii) unless $v$ has label $d$, and among these half-edge labels at most one of them is from the set $[m]$.
In particular, there is a fiber rooted on $d$ whose rightmost half-edge is labelled $1$.
Moreover, by construction, the last removed fiber must be attached to a child of root $d$.
Thus, label $(\overline{t+\ell_1^{\lambda}})^*$ must be the label for one of the children of $d$.
Regarding the sizes of this subset of fibers, there are two cases: (a) if the vertex $d$ has (out)degree one, then the size distribution 
is the multiset $\{1\} \bigcup \{\mu_i-1: \mu_i>1\}$, and (b) if the vertex $d$ has (out)degree $r> 1$, then the size distribution 
is the multiset $\{r\} \bigcup \bigl(\{\mu_i-1: \mu_i>1\} \setminus \{r-1\} \bigr)$. We shall refer to this subset of fibers as type-E fibers.
\end{itemize}
Next, we describe how to get back to a labelled plane tree from a set $F$ of fibers satisfying the above properties.
We view $F$ as a certain forest of trees.
\begin{itemize}
\item[(i)] Find a tree in $F$  with a minimum root such that there is no vertex with a starred label in the tree. If the root $v$ of the found tree has a label from $[d-1]$, then merge the root
with the vertex having the minimum label in the set $\{ (d+1)^*, (d+2)^*,\ldots, (d+\ell_1^{\mu})^*\}$ in $F$, and label the merged vertex by the label of $v$ (and discard the starred label before merging), let the newly attached half-edge of $v$ carry the `remembered' label;
If the root $v$ of the found tree has a label from $[\bar{t}]$, then merge the root
with the vertex having the minimum label in the set $\{ (\overline{t+1})^*, (\overline{t+2})^*,\ldots, (\overline{t+\ell_1^{\lambda}})^*\}$, and label the merged vertex by the label of $v$. Update $F$ as the resulting forest of trees.
\item[(ii)] iterate (i) until $F$ becomes a single labelled plane tree.
\end{itemize}
It should not be hard to verify the above constructions give a bijection between $\mathbb{T}$
and the set of forests under discussion.

Now we are ready to obtain $ |\mathbb{T}|$ by counting 
the number of distinct forests of fibers.
We may distinguish the following two cases:\\
\phantom\quad Case~$1$: the fiber rooted on $d$ has size $r=1$ provided that $\mu_i=1$ for some $i$.
First, we have ${t \choose \ell_1^{\lambda}}$ different ways to pick the labels of
the roots for the type-O fibers. We arrange these picked labels in increasing order.
Then, to guarantee the size distribution of these fibers, the sizes following the order will be 
a tuple $(t_1,\ldots, t_{\ell_1^{\lambda}})$ with its underlying supporting set being the multiset $\{\lambda_i-1: \lambda_i>1\}$. 
Next we determine the labels of the roots for the type-E fibers. We distinguish these roots into three classes:
the fiber with root $d$, the fibers where some of the half-edge labels in $[m]\setminus \{1\}$ appear, and the rest.
Suppose there are $b\leq \min\{m-1, \ell_1^{\mu}\}$ fibers in the second class.
Then there are ${d-1\choose b}$ ways to pick labels for the roots in the second class, and ${d-1-b \choose \ell_1^{\mu}-b}$
ways to pick labels for the roots of the third class.
Analogously, we arrange these roots of each class in increasing order,
and the sizes following the order will be the triple of tupes $\bigl((1), (r_1, \ldots, r_{b}), (q_1,\ldots, q_{\ell_1^{\mu}-b}) \bigr)$ such that
$$
  \{1\} \bigcup \{r_1, \ldots, r_{b}\} \bigcup \{q_1,\ldots, q_{\ell_1^{\mu}-b}\}=\{1\} \bigcup \{\mu_i-1 : \mu_i>1\}.
$$
Next we can arrange the remaining  $d-1+\ell_1^{\mu}-\ell_1^{\mu}$ unused labels from the set 
$[d-1]\bigcup 
\{ (d+1)^*, (d+2)^*,\ldots, (d+\ell_1^{\mu})^*\}$ to be the leaves in the type-O fibers.
Note that once the size distribution of these fibers are determined, different arrangements are just corresponding to different permutations.
So there gives $(d-1)!$ different ways. 
Next we need to assign half-edge labels to the half-edges incident to the leaves of the type-O fibers.
Note that the half-edge labels incident to the starred vertices are by construction binding to those `remembered'
ones. Thus we only need to assign half-edge labels to those unstarred vertices, among which $m-1-b$ of them are from $[m]\setminus \{1\}$.
There are ${m-1 \choose m-1-b} {n-m \choose (d-\ell_1^{\mu}-1)-(m-1-b)} (d-\ell_1^{\mu}-1)! $ different ways to do this.
For the unique fiber in the first class of the type-E fibers,
the unique child of the corresponding root must carry the label $(\overline{t+\ell_1^{\lambda}})^*$, and the label for the unique half-edge is by design $1$;
for the fibers in the second class of the type-E fibers,
we have to pick $\sum_j r_j$ labels from the unused labels in the set $[\bar{t}]\bigcup \{ (\overline{t+1})^*, (\overline{t+2})^*,\ldots, (\overline{t+\ell_1^{\lambda}})^*\}$, and arrange them linearly, which gives us ${t-1 \choose \sum_j r_j } (\sum_j r_j)!$ different ways.
We also need to assign half-edge labels from $[n]$ including the remembered ones.
This can be done by picking $\sum_j r_j$ unused labels from $[n]\setminus [m]$ first and arranging them linearly, which gives ${n-m- (d-\ell_1^{\mu}-1)+(m-1-b) \choose \sum_j r_j} (\sum_j r_j)!$
different ways, and next by inserting an unused element from $[m]\setminus \{1\}$ into each fiber for which we have $\prod_j (r_j+1) b!$ distinct ways.
For the third class, we just need to arrange the unused vertex labels linearly and the unused half-edge labels linearly, for which
there are $\bigl(\sum_j q_j \bigr)! \bigl(\sum_j (q_j+1) \bigr)!$ distinct ways. Thus, for this case, the total number of distinct forests of fibers is
\begin{multline*}
{t \choose \ell_1^{\lambda}}  \sum_{(t_1, \ldots, t_{\ell_1^{\lambda}})} (d-1)!  \sum_{b=0}^{\min\{m-1, \ell_1^{\mu}\}} {d-1\choose b} {d-1-b \choose \ell_1^{\mu}-b} 
\sum_{\substack{(r_1, \ldots, r_{b}),\\ (q_1,\ldots, q_{\ell_1^{\mu}-b})}} {m-1 \choose m-1-b} \\
\cdot {n-m \choose (d-\ell_1^{\mu}-1)-(m-1-b)} (d-\ell_1^{\mu}-1)!
{t-1 \choose \sum_j r_j } \left[\bigl(\sum_j r_j \bigr)! \right]^2  \prod_j (r_j+1) b! \\ 
\cdot {n-m- (d-\ell_1^{\mu}-1)+(m-1-b) \choose \sum_j r_j} \bigl(\sum_j q_j \bigr)! \bigl(\sum_j (q_j+1) \bigr)!.
\end{multline*}
\\
\phantom\quad Case~$2$: the fiber rooted on $d$ has size $r>1$ provided $\mu_i=r$ for some $i$.
First, we have ${t \choose \ell_1^{\lambda}}$ different ways to pick the labels of
the roots for the type-O fibers and we arrange these picked labels in increasing order. 
Next we determine the labels of the roots for the type-E fibers. We distinguish these roots into three classes:
the fiber with root $d$, the fibers where some of the half-edge labels in $[m]\setminus \{1\}$ appear, and the rest.
Suppose there are $b\leq \min\{m-1, \ell_1^{\mu}-1\}$ fibers in the second class.
Then there are ${d-1\choose b}$ ways to pick labels for the roots in the second class, and ${d-1-b \choose \ell_1^{\mu}-b-1}$
ways to pick labels for the roots of the third class.
Analogously, we arrange these roots of each class in increasing order,
and the sizes following the order will be the triple of tupes $\bigl((r), (r_1, \ldots, r_{b}), (q_1,\ldots, q_{\ell_1^{\mu}-b-1}) \bigr)$ such that
$$
  \{r\} \bigcup \{r_1, \ldots, r_{b}\} \bigcup \{q_1,\ldots, q_{\ell_1^{\mu}-b-1}\}=\{r\} \bigcup \bigl(\{\mu_i-1 : \mu_i>1\} \setminus \{r-1\} \bigr).
$$
Next we have $(d-1)!$ ways to assign labels for the leaves in the type-O fibers and
${m-1 \choose m-1-b} {n-m \choose (d-\ell_1^{\mu})-(m-1-b)} (d-\ell_1^{\mu})! $ different ways to 
assign half-edge labels to unstarred vertices there.
For the unique fiber in the first class of the type-E fibers,
besides the starred label $(\overline{t+\ell_1^{\lambda}})^*$ for one of the leaves,
we need to pick up $r-1$ other unused labels with bars and
arrange all these $r$ labels linearly, which gives us ${t-1 \choose r-1} r!$ options,
 and we also need to pick $r-1$ unused half-edge labels from $[n]\setminus [m]$ and arrange them linearly,
 which gives us ${n-m-(d-\ell_1^{\mu})+(m-1-b) \choose r-1} (r-1)!$ different possibilities;
for the fibers in the second class of the type-E fibers,
we have to pick $\sum_j r_j$ labels from the unused labels in the set $[\bar{t}]\bigcup \{ (\overline{t+1})^*, (\overline{t+2})^*,\ldots, (\overline{t+\ell_1^{\lambda}})^*\}$, and arrange them linearly, which gives us ${t-r \choose \sum_j r_j } (\sum_j r_j)!$ different ways.
We also need to assign half-edge labels including the remembered ones,
which gives ${n-m-r+1-(d-\ell_1^{\mu})+(m-1-b) \choose \sum_j r_j} (\sum_j r_j)! \prod_j (r_j+1)b!$ distinct ways.
For the third class, we have $\bigl(\sum_j q_j \bigr)! \bigl(\sum_j (q_j+1) \bigr)!$ distinct ways to assign labels to the leaves there. Thus, for this case, the total number of distinct forests of fibers is
\begin{multline*}
{t \choose \ell_1^{\lambda}}  \sum_{(t_1, \ldots, t_{\ell_1^{\lambda}})} (d-1)!  \sum_{b\geq 0}{d-1\choose b} {d-1-b \choose \ell_1^{\mu}-b-1} 
\sum_{\substack{r>1,\\(r_1, \ldots, r_{b}),\\ (q_1,\ldots, q_{\ell_1^{\mu}-b-1})}} {m-1 \choose m-1-b} (d-\ell_1^{\mu})! \\
\cdot {n-m \choose (d-\ell_1^{\mu})-(m-1-b)}  {t-1 \choose r-1}  r! 
{t-r \choose \sum_j r_j } \left[\bigl(\sum_j r_j \bigr)! \right]^2  \prod_j (r_j+1) b! (r-1)! \\ 
\cdot {n-(d-\ell_1^{\mu})-1-b \choose r-1} {n-r- (d-\ell_1^{\mu})-b\choose \sum_j r_j} \bigl(\sum_j q_j \bigr)! \bigl(\sum_j (q_j+1) \bigr)!  .
\end{multline*}
Summing up the two cases, we obtain $|\mathbb{T}|$.
Then, we have ${p}^{\lambda}_{m,k}=\frac{|\mathbb{T}|}{t! (d-1)!}$ which can be simplified into the expression in the theorem, completing the proof.
\qed

Let $\mathfrak{I}^{\lambda}_{m,k}$ be the set of plane permutations on $[n]$
where the diagonal has cycle-type $\lambda$ and the vertical has $k$ cycles with 
the elements in $[m]$ ($m\leq n$) being fixed points.
Denote ${I}^{\lambda}_{m,k}=|\mathfrak{I}^{\lambda}_{m,k}|$.
Reasoning analogously, we obtain the following 

\begin{theorem}\label{thm:gen2}
For $\lambda\vdash n$ and $n+1-\ell(\lambda)-k>0$, we have
	\begin{align}
	 	 I^{\lambda}_{m,k}=\frac{\sum_{j\geq 1} {k+2j-m\choose 2j+1} I^{\lambda}_{m,k+2j}+\sum_{j\geq 1} \sum_{\mu \rhd_{2j+1} \lambda} \kappa_{\mu,\lambda} I^{\mu}_{m,k}}{n+1-\ell(\lambda)-k} \; .
	\end{align}
\end{theorem}

 \begin{theorem}\label{thm:i-ini}
 Let 
 $\ell(\mu)=d$, $\ell(\lambda)=t$ where $d+t=n+1$, and let
 $\mu= 1^{b_1}2^{b_2} \cdots n^{b_n}, \; \ell_1^{\mu}= \sum_{i>1} b_i$ where $b_1\geq m$,
 and $\lambda=1^{a_1}2^{a_2} \cdots n^{a_n}$.
Then
   \begin{align}
{I}^{\lambda}_{m, \mu}=\frac{(t-1)! (d-1)! (n+1-m)! }{\prod_{i\geq 1} a_i !  \prod_{i>1} b_i! (d-m+1-\ell_1^{\mu})! }  \; .
 \end{align}

 \end{theorem}

Let $\sigma^{\lambda}_m$ (resp.~$\hat{\sigma}^{\lambda}_m$) denote the concerned separation (resp.~isolation) probability 
w.r.t.~the cycle-type $\lambda\vdash n$.
Obviously, we have
\begin{align*}
\sigma^{\lambda}_m =\frac{\sum_k p^{\lambda}_{m,k}}{(n-1)!(n-1)!}, \qquad \hat{\sigma}^{\lambda}_m =\frac{\sum_k I^{\lambda}_{m,k}}{(n-1)!(n-1)!}.
\end{align*}
Based on Theorem~\ref{thm:gen1} (resp.~Theorem~\ref{thm:gen2}) and their respective initial values, we can compute the separation (resp.~isolation) probability for any $\lambda$.

\section{Analogues of the Zagier-Stanley formula}\label{sec3}

Denote $C_m(n,k)$ the number of permutations on $[n]$ with $k$ cycles where the elements in $[m]$ are in distinct cycles.
The numbers $C_m(n,k)$ relate to their counterparts $C(n,k)$ as follows:
\begin{align}
C_m(n,k) =\sum_{d\geq 0}{d+m-1 \choose d}d! C(n-d-m, k-m).
\end{align}
Based on eq.~\eqref{eq:U}, we first have the following corollary.

\begin{corollary}\label{cor:exc}
The total number of exceedances of plane permutations in $\bigcup_{\lambda \vdash n} U^{\lambda}_{m,k}$
is given by
\begin{multline}
\sum_{a>0, \, \lambda \vdash n} a p^{\lambda, n, a}_{m,k} = (n-k) (n-1)! C_m(n,k) \\
-\sum_{j\geq 1} \left[ m {k+2j-m \choose 2j} + {k+2j-m \choose 2j+1}\right] (n-1)! C_m(n, k+2j).
\end{multline}
\end{corollary}
\proof It should not be hard to see that for a fixed upper horizontal,
the total number of plane permutations over all possible diagonal cycle-types and such that 
the vertical has $k$ cycles with the elements in $[m]$ in distinct cycles is exactly $C_m(n,k)$.
Then, summing over all $\lambda\vdash n$ on both sides of eq.~\eqref{eq:U} will complete the proof.\qed

On the other hand, we can directly easily count the total number of exceedances.
\begin{lemma}\label{lem:exc}
The total number of exceedances of plane permutations in $\bigcup_{\lambda \vdash n} U^{\lambda}_{m,k}$ is
\begin{align}
\sum_{a>0, \, \lambda\vdash n} a p^{\lambda, n, a}_{m,k} = (n-1)! \left[ {n-m \choose 2} + m (n-m) \right] C_m(n-1,k).
\end{align}
\end{lemma}
\proof 
Firstly, we observe that for any two fixed upper horizontals, the respective total number  
of exceedances of plane permutations with the two upper horizontals are equal (by relabelling argument). Thus, it suffices to count the total number of exceedances of plane permutations with the upper horizontal $(1\; 2\; \cdots \; n)$. The set of exceedances under study
 can be equivalently represented as the set consisting of pairs $\bigl(\pi, (x, \pi(x)) \bigr)$ where $\pi$ is a permutation of 
$k$ cycles on $[n]$ such that the elements in $[m]$ are contained in distinct cycles and $x < \pi(x)$. 
For a fixed pair $(x, y)$ such that $x<y$, there is obviously $C_m(n-1, k)$ permutations on $[n]$ such that $x$ is
an exeedance and $\pi(x)=y$. Note that among $x$ and $y$, at most one of them is contained in $[m]$.
Therefore, the number of exceedances for a fixed upper horizontal is $\left[ {n-m \choose 2} + m (n-m) \right] C_m(n-1,k)$ 
and the proof follows. \qed

Based on Corollary~\ref{cor:exc} and Lemma~\ref{lem:exc}, we have

\begin{proposition}\label{prop:recurrence} 
The following relation holds:
\begin{multline}
(n+1-k) \frac{2 (n-1)! C_m(n+1,k)}{(n+m)(n+1-m)} = (n-1)!  C_m(n,k)\\
+\sum_{j\geq 1} \left[ m {k+2j-m \choose 2j} + {k+2j-m \choose 2j+1}\right] \frac{2 (n-1)! C_m(n+1, k+2j)}{(n+m)(n+1-m)} \; .
\end{multline}
\end{proposition}

In order to proceed, we need the following corollary.

\begin{corollary}\label{cor:recurrence} 
Let $p^{(n)}_{m,k}=p^{\lambda}_{m,k}$ for $\lambda=n^1$. Then
		\begin{align}
	p_{m,k}^{(n)}=\frac{\sum_{j\geq 1} \left[m{k+2j-m\choose 2j}+{k+2j-m\choose 2j+1}\right] p_{m,k+2j}^{(n)}+(n-1)! C_m(n,k)}{n+1-k} \; ,
	\end{align}
	where $p^{(n)}_{m,n}=(n-1)!$.
\end{corollary}
\proof Note that for $\lambda=n^1$, any partition $\mu$ of $n$ satisfies $\mu\rhd_k \lambda$ for some $k$ and $\kappa_{\mu, \lambda}=1$.
Furthermore, for a fixed $k$, those $\mu$ such that $p^{\mu}_{m,k}\neq 0$ must have $\ell(\mu)$ of the same parity.
In particular, if $k$ is such that $p^{n^1}_{m,k}\neq 0$, then $\ell(\mu)=2j+1$ for some $j\geq 0$.
Thus
$$
\sum_{j\geq 0} \sum_{\mu \rhd_{2j+1} \lambda} \kappa_{\mu,\lambda} p^{\mu}_{m,k}=(n-1)! C_m(n,k).
$$
Based on Theorem~\ref{thm:gen1}, we obtain the formula for $p^{(n)}_{m,k}$. The formula for $p^{(n)}_{m,n}$ is clear, completing the proof. \qed

Now we are ready to present the following analogue of the Zagier-Stanley formula.
\begin{theorem}\label{thm:analog1}
 For $n-k$ being even, we have
\begin{align}\label{eq:analog1}
p^{(n)}_{m,k} =\frac{2 (n-1)! C_m(n+1,k)}{(n+m)(n+1-m)} \; .
\end{align}
\end{theorem}
\proof 
First we know that $p^{(n)}_{m,k}\neq 0$ iff $n-k$ is even.
Based on Proposition~\ref{prop:recurrence} and Corollary~\ref{cor:recurrence}, we observe that the quantities on both sides of eq.~\eqref{eq:analog1} satisfy the same recurrence.
Thus, it suffices to check the respective initial conditions, that is, for $k=n$. As to the left-hand side, we have $p^{(n)}_{m,n}=(n-1)!$.
On the other hand, we have
$$
C_m(n+1, n)={n+1-m\choose 2}+m(n+1-m)=\frac{(n+m)(n+1-m)}{2} \; ,
$$
i.e., besides $n-1$ fixed points, there is a cycle of length two such that at most one of the two elements is contained 
in $[m]$. So, the initial value for the right-hand side is also $(n-1)!$, whence eq.~\eqref{eq:analog1}. \qed

Denote by $\widehat{C}_m(n,k)$ the number of permutations on $[n]$ with $k$ cycles where the elements in $[m]$ are fixed points. Obviously, we have
$
\widehat{C}_m(n,k) =C(n-m, k-m).
$
With completely the same idea as above, we can obtain another analogue of the Zagier-Stanley result.

\begin{theorem}\label{thm:analog2} 
For $m<n$ and $n-k$ being even, we have
\begin{align}\label{eq:analog2}
I^{(n)}_{m,k} =\frac{2 (n-1)! \widehat{C}_m(n+1,k)}{(n-m)(n+1-m)} \; .
\end{align}
\end{theorem}

\section{Separation and isolation probabilities for $\lambda=n^1$}\label{sec4}

In this section, we derive the concerned separation and isolation probabilities for the case $\lambda=n^1$,
denoted by $\sigma^{(n)}_m$ and $\hat{\sigma}^{(n)}_m$, respectively,

\begin{theorem}
For $m \leq n$, the separation probability $\sigma^{(n)}_m$ is given by
\begin{align}
\sigma^{(n)}_m=
\begin{cases}
\frac{1}{m!}, & \mbox{if $n-m$ is odd,}\\
\frac{1}{m!}+\frac{2}{(m-2)! (n+1-m)(n+m)}, &  \mbox{if $n-m$ is even.}
\end{cases}
\end{align}
\end{theorem}
\proof For $n$ being odd, the number of cycles in the vertical must be odd. Thus, the separation probability
is given by
$$
\frac{\sum_k p^{(n)}_{m,k}}{(n-1)!(n-1)!}=\frac{\sum_{k\in odd}\frac{2 (n-1)! C_m(n+1,k)}{(n+m)(n+1-m)}}{(n-1)!(n-1)!}=\frac{2\sum_{k\in odd} C_m(n+1,k)}{(n+m)(n+1-m) (n-1)!}  \; .
$$
If $n-m$ is odd (i.e.,~$m$ even), additional odd number of cycles other than the $m$ cycles containing the elements in $[m]$ exist in the vertical. Then
\begin{align*}
{\sum_{k\in odd}C_m(n+1,k)}
&=\sum_{d=0}^{n-m} {n+1-m\choose d} d! {d+m-1\choose m-1} \sum_{k\in odd} C(n+1-m-d,k)\\
&=\frac{(n+1-m)!}{2} \sum_{d=0}^{n-m-1} {d+m-1\choose d}+ \frac{(n+1-m) (n-1)!}{(m-1)!}\\
&= \frac{(n+1-m)!}{2}  {n-1\choose n-m-1}+ \frac{(n+1-m) (n-1)!}{(m-1)!}\\
&=\frac{(n+m)(n+1-m) (n-1)!}{2 m!} \; ,
\end{align*}
which implies the separation probability to be $\frac{1}{m!}$.
Note that we used the fact that $\sum_{k\in odd} C(n,k)=\sum_{k\in even} C(n,k)=\frac{n!}{2}$ for $n\geq 2$.

Analogously, for $n-m$ being even, we have
\begin{align*}
{\sum_{k\in odd}C_m(n+1,k)}
&=\sum_{d=0}^{n+1-m} {n+1-m\choose d} d! {d+m-1\choose m-1} \sum_{k\in even} C(n+1-m-d,k)\\
&= \frac{(n+1-m)!}{2}  {n-1\choose n-m-1}+ \frac{n!}{(m-1)!}\\
&=\frac{(n+m)(n+1-m) (n-1)!}{2 m!}+\frac{(n-1)!}{(m-2)!} \; ,
\end{align*}
which implies the separation probability to be $\frac{1}{m!}+\frac{2}{(m-2)! (n+1-m)(n+m)}$.
The case of $n$ being even is analogous, completing the proof. \qed

\begin{theorem}\label{thm:iso}
For $m < n$, the isolation probability $\hat{\sigma}^{(n)}_m$ is given by
\begin{align}
\hat{\sigma}^{(n)}_m=\frac{1}{m!} {n-1 \choose m}^{-1}.
\end{align}
\end{theorem}
\proof For $n$ being odd, the number of cycles in the vertical must be odd. Thus, the isolation probability
is given by
$$
\frac{\sum_k I^{(n)}_{m,k}}{(n-1)!(n-1)!}=\frac{\sum_{k\in odd}\frac{2 (n-1)! \widehat{C}_m(n+1,k)}{(n-m)(n+1-m)}}{(n-1)!(n-1)!}=\frac{2\sum_{k\in odd} \widehat{C}_m(n+1,k)}{(n-m)(n+1-m) (n-1)!} \; .
$$
If~$m$ is even, then there must be an odd number of additional cycles other than the $m$ fixed points contained in $[m]$. Thus
\begin{align*}
{\sum_{k\in odd}\widehat{C}_m(n+1,k)}
= \sum_{k\in odd} C(n+1-m, k-m)=\frac{(n+1-m)!}{2} \; ,
\end{align*}
which implies the isolation probability to be $\frac{1}{m!} {n-1 \choose m}^{-1}$.
Analogous analysis on the remaining cases give the same isolation probability, whence the theorem.\qed

As consequences of Theorem~\ref{thm:iso}, we immediately obtain the following corollaries.

\begin{corollary}\label{fix-free}
The probability for the product of two uniformly random $n$-cycles being fixed point free
is given by $ \sum_{j\geq 0}^{n-1} (-1)^j \frac{n}{(n-j) j!} + (-1)^n \frac{1}{(n-1)!}$.
\end{corollary}
\proof Applying Theorem~\ref{thm:iso} and the inclusion-exclusion principle, we obtain that the total number of pairs
of $n$-cycles whose product is fixed point free is given by 
\begin{align*}
&(n-1)!(n-1)!  \left[ \sum_{j\geq 0}^{n-1} (-1)^j {n\choose j} \frac{1}{j!} {n-1 \choose j}^{-1} + (-1)^n \frac{1}{(n-1)!}\right]\\
= &(n-1)!(n-1)!  \left[ \sum_{j\geq 0}^{n-1} (-1)^j \frac{n}{(n-j) j!} + (-1)^n \frac{1}{(n-1)!}\right] \; ,
\end{align*}
whence the corollary. \qed

\begin{corollary}
The expected number of fixed points in the product of two uniformly random $n$-cycles is $\frac{n}{n-1}$,
and the variance of the number of fixed points is given by
$$
\sum_{i=0}^n i^2 \left( \sum_{j\geq i}^{n-1} (-1)^{j-i} \frac{n}{(n-j) j!} + (-1)^{n-i} \frac{1}{(n-1)!} \right)-\frac{n^2}{(n-1)^2} \; .
$$
\end{corollary}
\proof Let $s_1$ and $s_2$ be two uniformly random $n$-cycles.
For the random permutation $\pi=s_1 s_2$, for $1\leq i \leq n$, let $X_i$ be a random variable defined as follows: $X_i=1$ if
$i$ is a fixed point of $\pi$, $X_i=0$ otherwise. Based on Theorem~\ref{thm:iso}, we have the expectation of $X_i$, 
$\text{E}(X_i)=\frac{1}{n-1}$.
Note that the number of fixed points of $\pi$ is $X=\sum_i X_i$.
Therefore, the expected number of fixed points is given by
$$
\text{E} (X) =\sum_i \text{E}(X_i)=\frac{n}{n-1} \; ,
$$
completing the proof of the former part. As to the latter part,
applying Theorem~\ref{thm:iso} and the inclusion-exclusion principle, we obtain the total number of pairs
of $n$-cycles whose product has exactly $i$ fixed points.
Then, we can compute $\text{E}(X^2)$, and the variance follows
from computing $\text{E}(X^2) -\bigl( \text{E}(X) \bigr)^2$.
This completes the proof.
\qed


We remark that by
applying our result eq.~\eqref{eq:analog2}, we can also obtain the expectation and variance of fixed points
over products of 
pairs of long cycles whose product has exactly $k$ cycles.

We conclude this section by proving the mentioned elegant formula for another kind of separation probabilities.
Let $m>0$ and $\alpha=(\alpha_1, \alpha_2, \ldots, \alpha_m)$ be an integer composition of $n$, i.e.,~$\sum_i \alpha_i =n$ and $\alpha_i>0$. We write $\alpha \models n$.
Let $B_i=\{\sum_{j=0}^{i-1} \alpha_j+1, \ldots, \sum_{j=0}^{i-1} \alpha_j+\alpha_i\} \subseteq [n]$ where we assume $\alpha_0=0$.
A permutation $\pi$ on $[n]$ is called $\alpha$-separated if the elements in every cycle of $\pi$
are coming from the same $B_i$ for some $1\leq i \leq m$.

Based on Theorem~\ref{thm:iso}, for $\alpha=(1,\ldots, 1, n+1-m) \models n$, we obtain that the number of pairs of 
$n$-cycles whose product is $\alpha$-separated is ${(n-1)! (n-m)!}$.
In fact, we can extend this result to any $\alpha$. But we have to take advantage of the following 
symmetry property.

\begin{proposition}[Bernardi et al.~\cite{bdms}] \label{prop:symmetry}
Let $\alpha,~\beta \vdash n$ with $k$ components. 
Let $p^{(n)}_{\alpha}$ denote the number of pairs of 
$n$-cycles whose product is $\alpha$-separated.
Then, $\frac{p^{(n)}_{\alpha}}{p^{(n)}_{\beta}} =\frac{\prod_{i=1}^k \alpha_i !}{\prod_{i=1}^k \beta_i !}$.
\end{proposition}

Combined with our obtained formula for $\alpha=(1,\ldots, 1, n+1-m) \models n$, we obtain

\begin{proposition} For any $\alpha=(\alpha_1, \ldots, \alpha_k) \models n$, we have
$$
p^{(n)}_{\alpha}=\frac{\prod_{i=1}^k \alpha_i !}{(n+1-k)!} {(n-1)! (n-k)!} =\frac{ (n-1)!\prod_{i=1}^k \alpha_i !}{n+1-k} \; .
$$
\end{proposition}

\section{Final remarks}\label{sec6}

As pointed out by one anonymous referee,
another approach of computing the separation probability for a general $\lambda$
is as follows: for each cycle-type $\mu$, get the number of factorizations of a permutation of 
cycle-type $\mu$ into a long cycle and a permutation of cycle-type $\lambda$, e.g., by using the result
in~\cite{gs,cff,jackson}, and then multiply with the number of permutations of cycle-type $\mu$ separating the
elements in $[m]$, and finally sum over all cycle-types $\mu$.
However, none of these involved formulas are simple.
It is likely that
the overall computational complexity of our approach here is slightly smaller.

The anonymous referee also pointed out that the number of pairs of $n$-cycles whose product have $k$ cycles and such that the elements
in $[m]$ are fixed can be obtained via some modification of an argument of Cori, Marcus
and Schaeffer~\cite[Corollary~$1$]{cms}.

We remark that the result due to Du and Stanley as well as our Corollary~\ref{fix-free}
have been also obtained in B\'{o}na and Pittel~\cite{bona-pittel} by some computation
based on Fourier transform. F\'{e}ray and Rattan~\cite{fr} have also obtained
the formulas for a general $m$, mostly inductive.

\section*{Acknowledgments}

I thank one anonymous referee for some valuable comments and pointing me to the reference of B\'{o}na and Pittel.
I also thank Christian Reidys for encouragements and support.


\end{document}